\documentstyle[verbatim]{amsart}  

\newtheorem{theorem}{Theorem}[section]

\newtheorem{proposition}[theorem]{Proposition}
\newtheorem{corollary}[theorem]{Corollary}

\theoremstyle{definition}

\theoremstyle{remark}
\newtheorem{remark}[theorem]{Remark}

\numberwithin{equation}{section}

\begin{document}
\begin{flushright}
{\tt July 25, 2001\\
Revised September 28, 2001}
\end{flushright}

\title[Segal-Bargmann transforms, one-mode interacting Fock space]
{Segal-Bargmann Transforms of One-mode Interacting Fock Spaces
 Associated with \\ Gaussian and Poisson Measures}

\author[Nobuhiro Asai]{Nobuhiro Asai*}
\address{International Institute for Advanced Studies,
Kizu, Kyoto, 619-0225, Japan.}
\email{asai@@iias.or.jp, nobuhiro.asai@@nifty.com}
\thanks{*Research supported by a Postdoctoral Fellowship
of the International Institute for Advanced Studies, Kyoto, Japan}

\author{Izumi Kubo}
\address{Department of Mathematics,
Graduate School of Science,
Hiroshima University,
Higashi-Hiroshima, 739-8526, Japan}
\email{kubo@@math.sci.hiroshima-u.ac.jp}

\author{Hui-Hsiung Kuo}
\address{Department of Mathematics,
Louisiana State University,
Baton Rouge, LA 70803, USA.}
\email{kuo@@math.lsu.edu}

\date{July 25, 2001}

\begin{abstract}
Let $\mu_{g}$ and $\mu_{p}$ denote the Gaussian and Poisson 
measures on ${\Bbb R}$, respectively. We show that there
exists a unique measure $\widetilde{\mu}_{g}$ on ${\Bbb C}$
such that under the Segal-Bargmann transform $S_{\mu_g}$
the space $L^2({\Bbb R},\mu_g)$ is isomorphic to the space
${\cal H}L^2({\Bbb C}, \widetilde{\mu}_{g})$ of analytic 
$L^2$-functions on ${\Bbb C}$ with respect to 
$\widetilde{\mu}_{g}$. We also introduce the Segal-Bargmann 
transform $S_{\mu_p}$ for the Poisson measure $\mu_{p}$
and prove the corresponding result. As a consequence, when
$\mu_{g}$ and $\mu_{p}$ have the same variance, 
$L^2({\Bbb R},\mu_g)$ and $L^2({\Bbb R},\mu_p)$ are isomorphic 
to the same space ${\cal H}L^2({\Bbb C}, \widetilde{\mu}_{g})$ 
under the $S_{\mu_g}$ and $S_{\mu_p}$-transforms, respectively.  
However, we show that the multiplication operators by $x$ on
$L^2({\Bbb R}, \mu_g)$ and on $L^2({\Bbb R}, \mu_p)$ act 
quite differently on ${\cal H}L^2({\Bbb C}, \widetilde{\mu}_{g})$.
\end{abstract}

\maketitle


\section{Introduction}

Let $\mu$ be a probability measure on ${\Bbb R}$ having 
finite moments of all orders. In the paper \cite{ab}
Accardi and Bo\.zejko discovered a canonical
unitary isomorphism between the Hilbert space $L^2(\mu)$
and the one-mode interacting Fock space $\Gamma(\lambda)$
associated with a sequence $\lambda$ arising from $\mu$.
Under this isomorphism the number vectors $\Phi_n, n\geq 0,$
conrrespond to the orthogonal polynomials $P_n(x)$ 
associated with $\mu$ and the modified field operator
on $\Gamma(\lambda)$ corresponds to the multiplication 
operator by $x$ on $L^2(\mu)$.

Being motivated by Accardi-Bo\.zejko's discovery, 
Asai has recently introduced in \cite{asai01} coherent 
vectors which are used to define the $S_{\mu}$-transform. 
The $S_{\mu}$-transform 
is shown in \cite{asai01} to be a unitary operator from 
$L^2(\mu)$ onto a Hilbert space ${\cal H}_{\lambda}$ of analytic 
functions on a disk $\Omega_{\lambda}\subset {\Bbb C}$,
where $\lambda$ is determined by $\mu$. The composition
of the Accardi-Bo\.zejko isomorphism and the 
$S_{\mu}$-transform gives the interacting Fock space 
counterpart of the well-known Segal-Bargmann transform
(cf. \cite{barg1,barg2,kk,segal1,segal2}).

The purpose of this paper is to apply the results of
Accardi-Bo\.zejko \cite{ab} and Asai \cite{asai01} to
the cases of Gaussian measure $\mu_g$ and Poisson measure 
$\mu_p$.  In paticular, when $\mu_g$ and $\mu_p$  have 
the same variance, we will see that by
Theorems \ref{thm:gaussian} and \ref{thm:poisson}
the Segal-Bargmann
transforms $S_{\mu_g}$ and $S_{\mu_p}$ take $L^2(\mu_g)$ 
and $L^2(\mu_p)$, respectively, to the same space
${\cal H}L^2({\Bbb C},\widetilde{\mu}_{g})$ of analytic
$L^{2}$-functions, where $\widetilde{\mu}_{g}$ is the
Gaussian measure on ${\Bbb C}$. However, 
the Segal-Bargmann representation of 
multiplication by $x$ on $L^2(\mu_g)$ is quite different 
from that on $L^2(\mu_p)$.

\section{One-mode Interacting Fock space}\label{sec:ifs}

Let $\mu$ be a probability measure on ${\Bbb R}$ having 
finite moments of all orders. It is well-known 
\cite{szego} that there exist (1) a complete orthogonal 
system $\{P_{n}(x)\}_{n=0}^{\infty}$ of polynomials
for $L^{2}(\mu)$ with $P_{0}=1$, (2) a sequence 
$\{\omega_{n}\}_{n=1}^{\infty}$ of nonnegative real 
numbers, and (3) a sequence $\{\alpha_{n}\}_{n=0}^{\infty}$ 
of real numbers such that the following equalities hold 
for all $n\geq 0$:
\begin{align}
 (x-\alpha_{n})P_{n}(x)  & = P_{n+1}(x) + 
  \omega_{n} P_{n-1}(x), \label{eq:a} \\
  \langle P_{n}, P_{m}\rangle_{L^{2}(\mu)} 
   & = \delta_{n, m} \omega_{1} \cdots \omega_{n}, \notag
\end{align}
where $\omega_0P_{-1}=0$ by convention. We have the fact 
that the sequence $\alpha_{n}
=0$ for all $n$ if and only if $\mu$ is symmetric.

For a probability measure $\mu$ with the associated
sequence $\{\omega_{n}\}_{n=1}^{\infty}$, we define a
sequence $\{\lambda_{n}\}_{n=0}^{\infty}$ by
\begin{equation}\label{eq:lambda}
\lambda_0=1, \quad \lambda_{n}= 
\omega_{1}\omega_{2} \cdots \omega_{n}, \quad n\geq 1.
\end{equation}

Assume that the sequence $\{\lambda_{n}\}_{n=0}^{\infty}$ 
satisfies the condition:
\smallskip
\begin{equation} \label{eq:condition}
 \!\!\!  \!\!\! \!\!\! \!\!\!  \!\!\! \!\!\!
 \!\!\!  \!\!\! \!\!\! (\star) \qquad
 \inf_{n\geq 0} 
\lambda_{n}^{1\over n} > 0.  
\quad  \qquad \qquad\qquad\quad 
\end{equation}

\medskip
With such a sequence $\lambda=\{\lambda_{n}\}_{n=0}^{\infty}$, 
we define $\Gamma(\lambda)$ by
\begin{equation*}
 \Gamma(\lambda) = \Big\{(a_{0}, a_{1}, \ldots, a_{n},
  \ldots):\, a_{n}\in {\Bbb C}, \, \sum_{n=0}^{\infty}
  \lambda_{n} |a_{n}|^{2} < \infty \Big\}
\end{equation*}
and a norm $\|\cdot\|_{\lambda}$ on $\Gamma(\lambda)$ by
$$\|(a_{n})\|_{\lambda} = \bigg(\sum_{n=0}^{\infty}
 \lambda_{n} |a_{n}|^{2}\bigg)^{1/2}.
$$
Then $\Gamma(\lambda)$ is a Hilbert space with norm
$\|\cdot\|_{\lambda}$. It is called the
{\it one-mode interacting Fock space} associated with 
$\lambda$ \cite{ab,alv}.

Define a {\it number vector} $\Phi_n, n\geq 0$, by
$$
\Phi_n = (0, \ldots, 0, \overset{\overset{n+1}{\lor}}{1},
  0, \ldots).
$$

The vector $\Phi_{0}$ is called a {\it vacuum vector}.
Let $A$ be a densely defined operator on $\Gamma(\lambda)$
such that
$$
A\Phi_{0} = 0, \quad A\Phi_{n} = \omega_{n}\Phi_{n-1},
 \quad n\geq 1.
$$
The adjoint operator $A^{*}$ of $A$ is easily checked to be
given by
$$
A^{*}\Phi_{n} = \Phi_{n+1}, \quad n\geq 0.
$$
The operators $A$ and $A^{*}$ are called the 
{\it annihilation} and {\it creation operators} on
$\Gamma(\lambda)$, respectively. 
The number operator $N$ is defined by 
\begin{equation*}
N\Phi_n = n\Phi_n, \quad n\geq 0.
\end{equation*}

In addition, we define an 
operator $\alpha_N$ on $\Gamma(\lambda)$ by
\begin{equation*} \label{eq:lambdaop}
\alpha_N \Phi_{n} = \alpha_{n} \Phi_{n}, \quad
 n\geq 0.
\end{equation*}

Now we can state the result of Accardi and Bo\.zejko 
\cite{ab}: There exists an unitary isomorphism 
$U: \Gamma(\lambda) \to L^{2}(\mu)$ satisfying the
following conditions:

\begin{itemize}
\item[(1)] $U\Phi_{0} = 1$,
\item[(2)] $UA^{*}U^{*} P_{n} = P_{n+1}$, 
\item[(3)] $U(A+A^{*}+\alpha_N)U^{*} =Q$, where $Q$ is
the multiplication operator by $x$ on $L^2(\mu)$.
\end{itemize}

\section{Segal-Bargmann Transform} \label{sec:sg}

Let $\{\lambda_{n}\}_{n=0}^{\infty}$ be the sequence
defined in Equation (\ref{eq:lambda}) and consider the
following series of a complex number $z$:
\begin{equation} \label{eq:g}
 G_{\lambda}(z) = \sum_{n=0}^{\infty}\frac{1}
 {\lambda_n} z^n.
\end{equation}
Note that by Condition ($\star$) in Equation 
(\ref{eq:condition}) this series has a positive radius
of convergence, denoted by $r_{\lambda}$ . 

Let
$\Omega_{\lambda} = \{z\in {\Bbb C}: |z|< \sqrt{r_{\lambda}}\}$.
For each $z\in \Omega_{\lambda}$, Asai \cite{asai01} has
introduced a {\it coherent vector} 
$E_{\lambda}(\cdot, z)$ with respect to the family 
$\{P_{n}\}$ in Equation (\ref{eq:a}) by
\begin{equation} \label{eq:b}
 E_{\lambda}(x, z) = \sum_{n=0}^{\infty} {P_{n}(x) \over
 \lambda_{n}}\,z^{n}, \quad x\in {\Bbb R}.
\end{equation}

It is easy to see 
$$
\|E_{\lambda}(\cdot, z)\|_{L^{2}(\mu)} =
 G_{\lambda}(|z|^{2})^{1/2}
$$
and so $E_{\lambda}(\cdot, z) \in 
L^{2}(\mu)$ for all $z\in \Omega_{\lambda}$.
Moreover, the set $\{E_{\lambda}(\cdot, z): z \in 
\Omega_{\lambda}\}$ is linearly independent and spans a 
dense subspace of $L^{2}(\mu)$.

For $f\in L^{2}(\mu)$, let $S_{\mu} f$ be the function 
defined by:
\begin{equation} \label{eq:3-1}
(S_{\mu} f)(z)
 =\langle E_{\lambda}(\cdot, \overline{z}), 
 f \rangle_{L^2(\mu)} = \int_{{\Bbb R}} 
  E_{\lambda}(x, z) f(x)\,d\mu(x),
   \quad  z \in \Omega_{\lambda}.
\end{equation}

The mapping $S_{\mu}$ defined on $L^{2}(\mu)$ is called
the {\it Segal-Bargmann transform}. Asai has shown in
\cite{asai01} that $S_{\mu}$ is a unitary operator from
$L^{2}(\mu)$ onto ${\cal H}_{\lambda}$. Here 
${\cal H}_{\lambda}$ is given by
\begin{equation} \label{eq:hs}
 {\cal H}_{\lambda} = \bigg\{
 F(z) = \sum_{n=0}^{\infty} a_{n} z^{n}: \>
 F {\rm ~is~analytic~on~} \Omega_{\lambda}
 ~{\rm and~} \sum_{n=0}^{\infty} 
\lambda_{n} |a_{n}|^{2} < \infty\bigg\}.
\end{equation}
It is a Hilbert space with the norm
\begin{equation} \label{eq:hsnorm}
\|F\|_{\cal H_{\lambda}} = \left(\sum_{n=0}^{\infty} 
\lambda_{n} |a_{n}|^{2}\right)^{1/2}.
\end{equation}

Now, let us introduce operators $\widetilde{A}$ and 
$\widetilde{A}^{*}$ acting on ${\cal H}_{\lambda}$ by 
\begin{equation*}
\widetilde{A}1=0,
\quad \widetilde{A}z^n=\omega_nz^{n-1},\quad n\geq 1
\end{equation*}
and 
\begin{equation*}
\widetilde{A}^{*}z^n=z^{n+1},\quad n\geq 0.
\end{equation*}
Operators $\widetilde{A}$ and $\widetilde{A}^{*}$
satisfy the commutation relation 
$[\widetilde{A}, \widetilde{A}^{*}] z^{n} = 
(\omega_{n+1}-\omega_n) z^{n}$.  

The number operator $\widetilde{N}$
acting on ${\cal H}_{\lambda}$ is defined by 
\begin{equation*}
\widetilde{N}z^n=nz^n,\quad n\geq 0
\end{equation*}

In addition, we can define an operator $\widetilde{\alpha}_N$
acting on ${\cal H}_{\lambda}$ by 
\begin{equation*}
\widetilde{\alpha}_Nz^n=\alpha_nz^n,\quad n\geq 0 .
\end{equation*}
The operators 
$\widetilde{A}$, $\widetilde{A}^*$, $\widetilde{N}$ and 
$\widetilde{\alpha}_N$ correspond to
the operators 
$A$, $A^*$, $N$ and 
$\alpha_N$ on $\Gamma(\lambda)$, respectively.

\section{Main results} \label{sec:mr}

Our main results (Theorems \ref{thm:gaussian} and 
\ref{thm:poisson} below) in this papers are concerned 
with the special cases when $\mu$ is a Gaussian or 
Poisson measure.

\smallskip\noindent
{\bf Case 1}: $\mu=$ Gaussian measure $\mu_{g}$ with mean $m$ 
and variance $\sigma^2$.
\smallskip

From \cite{kuo96, szego} we have the following equalities:
\begin{align*}
 & H_{n+1}(x; \sigma^2) - x H_{n}(x; \sigma^2)
  + \sigma^2 n H_{n-1}(x; \sigma^2) = 0,  \label{eq:3-2} \\
 & \int_{{\Bbb R}} H_{n}(x; \sigma^2) H_{m}(x; \sigma^2)
 {1 \over \sqrt{2\pi \sigma^2}} e^{-x^{2}/2\sigma^2}\,dx
  = \delta_{n, m} \sigma^{2n} n! , \label{eq:3-3}
\end{align*}
where $H_{n}(x; \sigma^2)$ is the Hermite polynomial of
degree $n$ with parameter $\sigma^2$. Thus the three
quantities for the Gaussian measure $\mu_{g}$ in Equation
(\ref{eq:a}) are given by
\begin{align*}
 P_{n}(x) & = H_{n}(x-m; \sigma^2),  \\
 \omega_{n} & = \sigma^2 n,  \\
 \alpha_{n} & = m \quad \forall n \geq 1.
\end{align*}
Moreover, we have the associated quantities in Equations 
(\ref{eq:lambda}), (\ref{eq:g}), and (\ref{eq:b}):
\begin{align}
 \lambda_{n} & = \sigma^{2n} n!, \label{eq:3-4}  \\
 G_{\lambda}(z) & = \exp\Big[{z \over \sigma^2}\Big],
  \quad r_{\lambda} = \infty,  \label{eq:3-5} \\
 E_{\lambda}(x, z) & = 
 \exp\Big[{z \over \sigma^2} 
   (x-m) - {z^{2} \over 2\sigma^2}\Big]. \notag  
\end{align}
Obviously, Condition ($\star$) in Equation 
(\ref{eq:condition}) is satisfied.

\smallskip\noindent
{\bf Case 2}: $\mu=$ Poisson measure $\mu_{p}$ with  
parameter $a$.
\smallskip

From \cite{chihara} we have the following equalities:
\begin{align*}
 & C_{n+1}(x; a) = (x-n-a) C_{n}(x; a)
  - a n C_{n-1}(x; a),  \\
 & \int_{{\Bbb R}} C_{n}(x; a) C_{m}(x; a)\,
   d\mu_{p}(x) = \delta_{n, m}  a^{n} n!,
\end{align*}
where $C_{n}(x; a)$ is the Charlier polynomial of
degree $n$ with parameter $a$. Thus the three
quantities for the Poisson measure $\mu_{p}$ in Equation
(\ref{eq:a}) are given by
\begin{align}
 P_{n}(x) & = C_{n}(x; a), \notag \\
 \omega_{n} & = a n,  \notag \\
 \alpha_{n} & = n+a.  \notag
\end{align}
Moreover, we have the associated quantities in Equations 
(\ref{eq:lambda}), (\ref{eq:g}), and (\ref{eq:b}):
\begin{align*}
 \lambda_{n} & = a^{n} n!,  \\
 G_{\lambda}(z) & = \exp\Big[{z \over a}\Big],
  \quad r_{\lambda} = \infty,  \\
 E_{\lambda}(x, z) & = e^{-z}\Big(1+{z\over a}
        \Big)^{x}. 
\end{align*}
Obviously, Condition ($\star$) in Equation 
(\ref{eq:condition}) is satisfied.

For the Poisson case, we have 
\begin{equation}\label{eq:alpha}
 \alpha_N = {1\over a}A^{*}A + a
 = N+a,
\end{equation}
which implies that 
\begin{equation}\label{eq:p-decom}
 A^{*}+A+ \alpha_N =
 \biggl({1\over\sqrt{a}}A^{*}+\sqrt{a}\biggr)
 \biggl({1\over\sqrt{a}}A+\sqrt{a}\biggr).
\end{equation}

\begin{proposition} \label{pro:gaussian}
For the Gaussian measure $\mu_{g}$ on ${\Bbb R}$
with mean $m$ and variance $\sigma^2$ , the following 
equalities hold:
\begin{itemize}
\item[(a)] $S_{\mu_g}H_n (\cdot - m; \sigma^2) = z^{n}$,
\item[(b)] $S_{\mu_g}UAU^{*}H_n (\cdot - m; \sigma^2) 
  = \sigma^2nz^{n-1}$,
\item[(c)] $S_{\mu_g}UA^{*}U^{*}H_n (\cdot - m; \sigma^2)
  =z^{n+1}$,
\item[(d)] $S_{\mu_g}\big((x-m) H_n (\cdot - m; \sigma^2)\big)
  = z^{n+1} + \sigma^2nz^{n-1}$.
\end{itemize}
\end{proposition}

\begin{pf}
Conclusion (a) follows from Equations (\ref{eq:b}),
(\ref{eq:3-1}), (\ref{eq:3-3}), and (\ref{eq:3-4}). 
(b) follows from (a) and the fact that 
$UAU^{*} H_{n}(\cdot - m; \sigma^2) = 
\sigma^2n H_{n-1}(\cdot - m; \sigma^2)$.
(c) follows from (a) and the fact that
$UA^{*}U^{*} H_{n}(\cdot - m; \sigma^2) = H_{n+1}(\cdot - m; \sigma^2)$.
(d) follows from (a) and Equation (\ref{eq:3-2}).
\end{pf}

\begin{proposition} \label{pro:poisson}
For the Poisson measure $\mu_{p}$ on ${\Bbb R}$
with parameter $a$, the following 
equalities hold:
\begin{itemize}
\item[(a)] $S_{\mu_p}C_n (\cdot; a) = z^{n}$,
\item[(b)] $S_{\mu_p}UAU^{*}C_n (\cdot; a) 
  = anz^{n-1}$,
\item[(c)] $S_{\mu_p}UA^{*}U^{*}C_n (\cdot; a)
  =z^{n+1}$,
\item[(d)] $S_{\mu_p}\big(x C_n (\cdot; a)\big)
  = z^{n+1} + (n+a) z^{n} + anz^{n-1}$.
\end{itemize}
\end{proposition}

\begin{pf}
The idea is similar to the proof of Proposition 
\ref{pro:gaussian} with the Hermite polynomials 
being replaced by the Charlier polynomials.
\end{pf}

We point out that for the Poisson case, we have the equalities
for operators acting on ${\cal H}_{\lambda}$

\begin{equation*}
 \widetilde{\alpha}_N = 
 {1\over a}\widetilde{A}^{*}\widetilde{A} + a
 =\widetilde{N}+a
\end{equation*}
and 
\begin{equation*}
 \widetilde{A}^* + \widetilde{A} + \widetilde{\alpha}_N =
 \biggl({1\over\sqrt{a}}\widetilde{A}^{*}+\sqrt{a}\biggr)
 \biggl({1\over\sqrt{a}}\widetilde{A}+\sqrt{a}\biggr)
\end{equation*}
which correspond to the Equations \eqref{eq:alpha} and 
\eqref{eq:p-decom}, respectively.

Then we can apply Propositions 
\ref{pro:gaussian} and \ref{pro:poisson} to get the 
next result. 
\begin{corollary}
Let $\mu_{g}$ be the Gaussian measure on 
${\Bbb R}$ with 
mean $m$ and variance $\sigma^2$. 
Let $\mu_p$ be the Poisson measure on 
${\Bbb R}$ with the parameter $a$.
Then
\begin{itemize}
\item[(1)] $S_{\mu_g}\big((x-m) H_n (\cdot - m; \sigma^2)\big)
  = (\widetilde{A}^*+\widetilde{A})
  S_{\mu_g}H_n (\cdot - m; \sigma^2)$.
\item[(2)] $S_{\mu_p}\big(x C_n (\cdot; a)\big)
  = (\widetilde{A}^*+\widetilde{A}+
  \widetilde{\alpha}_N)S_{\mu_p}C_n (\cdot; a)$.
\end{itemize}
\end{corollary}

Now, suppose $\mu$ is a probability measure on ${\Bbb R}$
with finite moments of all orders. Let $\lambda$ be the
sequence associated with $\mu$ as given in Equation
(\ref{eq:lambda}). With this $\lambda$, we have a Hilbert
space ${\cal H}_{\lambda}$ of analytic functions on
$\Omega_{\lambda}$ in Equation (\ref{eq:hs}) with
norm $\|\cdot\|_{{\cal H}_{\lambda}}$ in Equation
(\ref{eq:hsnorm}). Consider the following 

\smallskip\noindent
{\bf Question:} Does there exist a unique measure 
$\widetilde\mu$ on $\Omega_{\lambda}$ such that 
$F\in {\cal H}_{\lambda}$ if and only if $F$ is analytic 
on $\Omega_{\lambda}$ and $F\in L^{2}(\Omega_{\lambda}, 
\widetilde\mu)$ with
\begin{equation} \label{eq:4-1}
 \|F\|_{{\cal H}_{\lambda}}^{2} =
  \int_{\Omega_{\lambda}} |F(z)|^{2}\,d\widetilde\mu(z) \quad ?
\end{equation}

\smallskip
Bargmann \cite{barg1} considered the equality
in Equation \eqref{eq:4-1} for the multidimensional 
standard Gaussian case. 
See also the paper by Gross and Malliavin \cite{gm}.

Our main results answer the above 
question for Gaussian and
Poisson measures.   
For convenience, let ${\cal H}L^{2}(\Omega_{\lambda}, 
\widetilde\mu)$ denote the Hilbert space of analytic 
functions $F$ on $\Omega_{\lambda}$ which are square 
integrable with respect to $\widetilde\mu$. 
The norm on
${\cal H}L^{2}(\Omega_{\lambda}, \widetilde\mu)$ is the
$L^{2}(\widetilde\mu)$-norm.

To answer the above question, 
we consider a criterion to check whether a measure $\nu$ 
satsifying the equation 
\begin{equation}\label{eq:criterion}
\int_{\Bbb C}\overline{z}^mz^nd\nu(z)=\gamma_{m,n}.
\end{equation}
is unique for given moments $\{\gamma_{m,n}\}$.

\begin{proposition}\label{prop:criterion}
Suppose $\{\gamma_{m,n}\}$ satsifies the condition
\begin{equation}\label{eq:gamma-condition}
\lim_{n\to \infty}{\gamma_{n,n}^{1/n}\over n^2}=0.
\end{equation}
Then the measure $\nu$ satisfying Equation \eqref{eq:criterion}
is unique.
\end{proposition}

\begin{pf}
Apply the Schwartz inequality to get
\begin{equation*}
|\gamma_{m,n}|=\biggl|\int_{\Bbb C}
\overline{z}^mz^nd\nu(z)\biggr|
\leq (\gamma_{m,m}\gamma_{n,n})^{1\over 2}.
\end{equation*}
Therefore the function 
\begin{equation*}
g(t,s)=\sum_{m,n=0}^{\infty}\gamma_{m,n}
{t^ms^n\over m!n!}
\end{equation*}
converges absolutely for $t,s\in {\Bbb C}$ by Equation 
\eqref{eq:gamma-condition}.
In fact, we see that 
\begin{align}\label{eq:4-11}
\int_{\Bbb C}\sum_{m,n=0}^{\infty}
{|(t\overline{z})^m(sz)^n|\over m!n!}d\nu(z)
& \leq 
\sum_{m,n=0}^{\infty}
{(\gamma_{m,m}\gamma_{n,n})^{1\over 2}\over m!n!}|t|^m|s|^n \notag\\
& \leq \biggl(\sum_{n=0}^{\infty}
{\gamma_{n,n}^{1\over 2}\over n!}R^n\biggr)^2
<\infty
\end{align} 
for $t, s\in {\Bbb C}$, $|t|,|s|\leq R$.  
Equation \eqref{eq:4-11} implies that 
$\exp[t\overline{z}+sz]$ is integrable with respect to a 
measure $\nu$ and 
\begin{equation*}
\int_{\Bbb C}\exp[{t\overline{z}+sz}]d\nu(z)=g(t,s)
\end{equation*} 
holds.  Therefore the characteristic function of a measure
$d\nu(x,y)=d\nu(z)$, $z=x+iy$, satisfies the equality 
\begin{equation*}
\int_{{\Bbb R}^2}\exp[i\xi x+i\eta y]d\nu(x,y)
=g\bigg({i\xi+\eta\over 2},{i\xi-\eta\over 2}\bigg)
\end{equation*}
for any $\xi, \eta\in {\Bbb R}$.  Hence $\nu$ is unique.
\end{pf}

\begin{remark} 
We are unable to find any literature dealing with the
moment problem for measures on the complex plane ${\Bbb C}$.
The above proposition gives a sufficient condition for the 
uniqueness of a measure on ${\Bbb C}$ in the moment problem.
On the other hand, many authors have studied 
the moment problem for measures on 
${\Bbb R}$. 
\end{remark}

\begin{theorem} \label{thm:gaussian}
Let $\mu_{g}$ be the Gaussian measure with mean $m$ and
variance $\sigma^2$. Let ${\cal H}_{\lambda}$ be the Hilbert 
space associated with $\mu_{g}$ as in Equation
(\ref{eq:hs}), i.e., $\Omega_{\lambda} = {\Bbb C}$ and
$\lambda_{n} = \sigma^{2n} n!$. Then there exists a unique
measure $\widetilde\mu$ on ${\Bbb C}$ such that
${\cal H}_{\lambda} = {\cal H} L^{2}({\Bbb C},
\widetilde\mu)$ and Equation (\ref{eq:4-1}) holds.
\end{theorem}

\begin{pf}
Consider the above question for the Gaussian measure
$\mu_{g}$ with mean $m$ and variance $\sigma^2$. In this case, 
we have $\Omega_{\lambda} = {\Bbb C}$ by Equation
(\ref{eq:3-5}) and $\lambda_{n} = \sigma^{2n}n!$ by Equation
(\ref{eq:3-4}). 

The uniqueness of $\widetilde{\mu}$ follows from 
Proposition \ref{prop:criterion}. Thus we only need to find
a measure $\widetilde\mu$ on ${\Bbb C}$ such that 
\begin{equation} \label{eq:4-2}
 \int_{{\Bbb C}}\overline{z}^mz^n\,d\widetilde{\mu}(z)
 =\delta_{m,n} \sigma^{2n}n!.
\end{equation}
Suppose $\widetilde{\mu}$ is given by $d\widetilde{\mu} 
(z) = {1\over 2\pi} d\rho(r)d\theta$ for $z=re^{i\theta}$. 
Then Equation (\ref{eq:4-2}) becomes
\begin{equation*}  
 \frac{1}{2\pi}\int_{0}^{\infty}
 \biggl(\int_{0}^{2\pi}
 e^{-i(m-n)\theta}\,d\theta\biggr)
 r^{m+n}\,d\rho(r)
 =\delta_{m,n} \sigma^{2n}n!.
\end{equation*}
This equality is obviously valid when $m\ne n$. Thus we are
looking for $\rho$ such that
\begin{equation}\label{eq:s-moment}
\int_{0}^{\infty} r^{2n}\,d\rho(r) = \sigma^{2n} n!.
\end{equation}
But it is easy to see that $\rho$ is given by
\begin{equation} \label{eq:cgm}
d\rho(r) = {2\over \sigma^2}\,
 r \exp\Big[{-{r^2\over \sigma^2}}\Big]\,dr.
\end{equation}

\noindent
The resulting measure 
$$
d\widetilde\mu_{g} (z) = 
{1\over \pi\sigma^2}\,r \exp\left[{-{r^2\over \sigma^2}}\right]
\,drd\theta \quad z = re^{i\theta}
$$ 
is a {\it Gaussian measure} on ${\Bbb C}$. 
Hence we have proved the existence of a measure $\widetilde\mu$ 
satisfying Equation \eqref{eq:4-1}.
\end{pf}

\begin{remark}\label{rem:4-7}
The measure given in Equation \eqref{eq:cgm} is the unique 
measure satisfying Equation \eqref{eq:s-moment}
since the sequence $\lambda_{n} = \sigma^{2n} n!$
satisfies the condition 
\begin{equation*}
\sum_{n=1}^{\infty}
(\lambda_n)^{-{1\over 2n}} = \infty 
\end{equation*}
in Theorem 1.11 of the book by Shohat and Tamarkin \cite{st}. 
However, in general it is not true that a measure 
$d\widetilde{\mu}(z)$ on ${\Bbb C}$ can be written as  
${1\over 2\pi} d\rho(r)d\theta$.  Thus we really need
Proposition \ref{prop:criterion} to show the uniqueness of
$\widetilde{\mu}$.
\end{remark}

Observe that we can also apply the above arguments to
the Poisson measure $\mu_{p}$. Thus we have the 
corresponding theorem for $\mu_{p}$.

\begin{theorem} \label{thm:poisson}
Let $\mu_{p}$ be the Poisson measure with parameter
$a$. Let ${\cal H}_{\lambda}$ be the Hilbert 
space associated with $\mu_{p}$ as in Equation
(\ref{eq:hs}), i.e., $\Omega_{\lambda} = {\Bbb C}$ and
$\lambda_{n} = a^{n} n!$. Then there exists a unique
measure $\widetilde\mu$ on ${\Bbb C}$ such that
${\cal H}_{\lambda} = {\cal H} L^{2}({\Bbb C},
\widetilde\mu)$ and Equation (\ref{eq:4-1}) holds.
\end{theorem}

In fact, it is easy to see that the unique measure 
$\widetilde\mu$ for $\mu_{p}$ is also a Gaussian measure
on ${\Bbb C}$ with $a$ replacing $\sigma^2$ in Equation
(\ref{eq:cgm}). In particular, we see that when $a=\sigma^2$
(i.e., $\mu_{g}$ and $\mu_{p}$ have the same variance)
the Segal-Bargmann representing spaces for $\mu_{g}$
and $\mu_{p}$ are the same space, namely, the space
${\cal H} L^{2}({\Bbb C}, \widetilde\mu_{g})$. However,
the multiplication operators for $\mu_{g}$ and $\mu_{p}$
are decomposed quite differently 
on ${\cal H} L^{2}({\Bbb C}, \widetilde\mu_{g})$.  

\begin{remark}
Hudson-Parthasarathy \cite{hp} and 
Ito-Kubo \cite{ik} obtained similar results presented in 
this section from the viewpoint of quantum stochastic calculus and 
white noise calculus (cf. \cite{kuo96}), respectively.
On the relationship between 
infinite dimensional Gaussian analysis 
and the Segal-Bargmann transform,
see the papers by Asai-Kubo-Kuo \cite{akk1,akk3,akk6,akk7}, 
Cochran-Kuo-Sengupta \cite{cks},
Kubo-Yokoi \cite{ky}, Lee \cite{lee}, Segal \cite{segal1, segal2},
Yokoi \cite{yokoi95} and the references therein.
\end{remark}

\begin{remark}
Recently, the $q$-deformed versions of Theorems \ref{thm:gaussian} and 
\ref{thm:poisson} have been
discussed by the first author \cite{asai02} following 
the technique in \cite{asai01}.
We remark that the case of $q$-deformed Gaussian measure 
has been investigated earlier by van Leeuwen-Maassen \cite{lm95}.
However, their method of analysis is slightly different from that in 
\cite{asai01, asai02}.  In addition, Kr\.olak 
examined the Segal-Bargmann space associated with $q$-commutation 
relation for $q>1$.
\end{remark}

\section*{Acknowledgments}
The authors would like to thank the referee 
for pointing out several relevant papers.
The results in this paper were obtained during several visits
of Asai and Kuo to Hiroshima University. Asai and Kuo
want to give their deep appreciation to Professor I. Kubo
for his invitation and kindness during their visits. They are 
thankful to the Mathematics Department of Hiroshima
University for the warm hospitality. 
Kuo is grateful for the financial support of 
Monbu-Kagaku-Sho (Ministry of Education and Science)
during his visit May 20--Aug 19, 2001.
\bibliographystyle{amsplain}

\begin{thebibliography}{99}

\bibitem{ab}L. Accardi and M. Bo\.{z}ejko,
\textit{Interacting Fock space and Gaussianization
of probability measures.}
Infinite Dimensional Analysis, Quantum Probability
and Related Topics,
\textbf{1}
(1998), 663--670.

\bibitem{alv} L. Accardi, Y.-G. Lu, and I. Volovich,
\textit{The QED Hilbert module and interacting Fock spaces.}
IIAS reports
\textbf{1997-008},
Pub. of IIAS (Kyoto),
1997.

\bibitem{asai01} N. Asai,
\textit{Analytic characterization of one-mode interacting
Fock space.}
Infinite Dimensional Analysis, Quantum Probability
and Related Topics, 
\textbf{4}
(2001), 409--415.

\bibitem{asai02} N. Asai,
\textit{Integral transform and Segal-Bargmann
representation associated to q-Charlier polynomials.}
{\em Preprint} (2001),
{\tt http://arXiv.org/abs/math.CA/0104260}

\bibitem{akk1} N. Asai, I. Kubo, and H.-H. Kuo,
\textit{Bell numbers, log-concavity, and log-convexity.}
Acta Appl. Math., 
\textbf{63} (2000), 79--87.

\bibitem{akk3} N. Asai, I. Kubo, and H.-H. Kuo,
\textit{CKS-space in terms of growth functions.}
in: {Quantum Information II},
T. Hida and K. Sait\^o ~(eds.)
World Scientific, 2000, pp. 17--27.

\bibitem{akk6} N. Asai, I. Kubo, and H.-H. Kuo,
\textit
{Roles of log-concavity, log-convexity, and growth order 
in white noise analysis.} Infinite Dimensional Analysis, 
Quantum Probability, and Related Topics,
\textbf{4} 
(2001) 59--84.

\bibitem{akk7} N. Asai, I. Kubo, and H.-H. Kuo,
\textit{General characterization theorems and intrinsic topologies
in white noise analysis.}
Hiroshima Math. J., 
\textbf{31} (2001), 
299--330.

\bibitem{barg1} V. Bargmann,
\textit{On a Hilbert space of analytic functions and
an associated integral transform, I.}
Comm. Pure Appl. Math.,
\textbf{14}
(1961), 187--214.

\bibitem{barg2} V. Bargmann,
\textit{On a Hilbert space of analytic functions and
an associated integral transform, II.}
Comm. Pure Appl. Math.,
\textbf{20}
(1967), 1--101.

\bibitem{chihara} T. S. Chihara,
\textit{An Introduction to Orthogonal Polynomials.}
Gordon and Breach, 1978.

\bibitem{cks} W. G. Cochran, H.-H. Kuo, and A. Sengupta,
\textit{A new class of white noise generalized functions.}
Infinite Dimensional Analysis, Quantum Probability
and Related Topics,
\textbf{1}
(1998), 43--67.

\bibitem{gm} L. Gross and P. Malliavin,
\textit{Hall's transform and the Segal-Bargmann map.}
in: It\^o Stochastic Calculus and Probability Theory,
N. Ikeda et al.~(eds.) Springer-Verlag, 1996,
pp. 73--116

\bibitem{hp} R. L. Hudson and K. R. Parthasarathy,
\textit{Quantum Ito's formula and stochastic evolutions.}
Comm. Math. Phys.,
\textbf{93} 
(1984), 301--323.

\bibitem{ik} Y. Ito and I. Kubo,
\textit{Calculus on Gaussian and Poisson white noises.}
Nagoya Math. J.,
\textbf{111}
(1988), 41--84.

\bibitem{Krolak} I. Kr\.olak,
\textit{Measures connected with Bargmann's representation of 
the $q$-commutation relation for $q>1$.}
Banach Center Publ.,
\textbf{43}
(1998), 253--257.

\bibitem{kk} I. Kubo and H.-H. Kuo,
\textit{Finite dimensional Hida distributions.}
J.~Funct.~Anal.,
\textbf{128}
(1995), 1--47.

\bibitem{ky} I. Kubo and Y. Yokoi,
\textit{Generalized functions and fluctuations in
fluctuation analysis.}
in: {Mathematical Approach to Fluctuations, Vol.II},
T. Hida et al. (eds.)
World Scientific, 1993, pp. 203--230.

\bibitem{kuo96} H.-H. Kuo,
\textit
{White Noise Distribution Theory.}
CRC Press, 1996.

\bibitem{lee} Y.-J. Lee,
\textit{Analytic version of test functionals, Fourier transform
and a characterization of measures in white noise calculus.}
J.~Funct.~Anal.,
\textbf{100}
(1991), 359--380.

\bibitem{lm95} H. van Leeuwen and H. Maassen, 
\textit{A $q$ deformation of the Gauss distribution.}
J. Math. Phys.,
\textbf{36}
(1995), 4743--4756.

\bibitem{segal1} I. E. Segal,
\textit
{Mathematical characterization of the physical vacuum for
a linear Bose-Einstein field.}
Illinois J. Math.,
\textbf{6}
(1962), 500--523.

\bibitem{segal2} I. E. Segal,
\textit
{The complex wave representation of the free Boson field.}
in: Essays Dedicated to M. G. Krein on
the Occassion of His 70th Birthday,
Advances in Math.: Supplementary Studies Vol.3,
I. Goldberg and M. Kac (eds.)
Academic, 1978, pp. 321--344.

\bibitem{st} J. Shohat and J. Tamarkin,
\textit{The Problem of Moments}.
Math Surveys \textbf{1},
Amer. Math. Soc., 1943.

\bibitem{szego} M. Szeg\"{o},
\textit{Orthogonal Polynomials}.
Coll. Publ. \textbf{23},
Amer. Math. Soc., 1975.

\bibitem{yokoi95} Y. Yokoi,
\textit{Simple setting for white noise calculus using
Bargmann space and Gauss transform.}
Hiroshima Math. J.,
\textbf{25}
(1995), 97--121.

\end{thebibliography}

\end{document}